\documentclass[letterpaper, 10 pt, conference]{ieeeconf}  

\IEEEoverridecommandlockouts                              
\usepackage{graphicx}
\usepackage{amsmath}
\usepackage{amssymb}

\newcommand{\diag}{\mathrm{diag}}
\newcommand{\hx}{{\hat{x}}}

\newcommand{\bDel}{{\bar\Delta}}
\newcommand{\B}[1]{{\bf #1}}
\newcommand{\Sc}[1]{{\mathcal{#1}}}
\newcommand{\R}[1]{{\rm #1}}
\newcommand{\mB}[1]{{\mathbb{#1}}}

\newcommand{\set}[2]{\left\{#1\,\left\vert\, #2\right.\right\}}

\newcommand{\half}{\frac{1}{2}}
\newcommand{\map}[3]{#1\,:\,#2\rightarrow #3}
\newcommand{\cS}{\mathcal{S}}
\newcommand{\cL}{\mathcal{L}}
\newcommand{\sfrac}[2]{\mbox{\small $\frac{#1}{#2}$}}

\usepackage{color}
\definecolor{softblue}{rgb}{0.90,0.92,1.00}

\newtheorem{lem}{Lemma}[section]
\newtheorem{remark}[lem]{Remark}

\newtheorem{algor}[lem]{Algorithm}

\newtheorem{lemma}{Lemma}[section]
\newtheorem{theorem}[lemma]{Theorem}
\newtheorem{corollary}[lemma]{Corollary}

\title{\LARGE \bf
Smoothing Dynamic Systems with State-Dependent Covariance Matrices
}

\author{Aleksandr Y. Aravkin and James V. Burke
\thanks{A. Y. Aravkin is with  IBM T.J. Watson Research Center, Yorktown Heights, 10598, NY, USA  
{\tt\small saravkin@us.ibm.com}}%
\thanks{J. V. Burke is with the Department of Mathematics, University of Washington, 
Seattle, WA, USA
{\tt\small jvburke@uw.edu}}%
}

\begin{document}

\maketitle

\begin{abstract}
Kalman filtering and smoothing algorithms are used in many areas, including 
tracking and navigation, medical applications, and financial trend filtering. 
One of the basic assumptions required to apply the Kalman smoothing framework
is that error covariance matrices are known and given. 
In this paper, we 
study a general class of inference problems where covariance matrices 
can depend functionally on unknown parameters. 
In the Kalman framework, this allows modeling situations 
where covariance matrices may depend functionally 
on the state sequence being estimated. We present an extended formulation 
and generalized Gauss-Newton (GGN) algorithm for inference in this context. 
When applied to dynamic systems inference, we show the algorithm  
can be implemented to preserve the computational efficiency of the classic Kalman smoother. 
The new approach is illustrated with a synthetic numerical example.

\end{abstract}

\section{Introduction}
\label{EGintroduction}
The Kalman filter~\cite{kalman} and smoother~\cite{RTS} are efficient algorithms 
to estimate the state of a dynamic system given noisy measurements. 
Over the last 10 years, the optimization perspective on the 
smoothing problem has produced many extensions to dynamic system estimation, 
including methods for smoothing systems with nonlinear process and 
measurement models~\cite{Bell1994}, systems with nonlinear inequality
constraints~\cite{Bell2008}, robust Kalman smoothing~\cite{RobustSparseArxiv, StArxiv, Aravkin2011tac},
and smoothing of sparse systems~\cite{Angelosante2009}.

In all of the above extensions, the variances of process and measurement errors
are assumed to be fixed and known. In practice, these quantities are often not known, 
and may in fact depend on the state. 
For example, radar position errors are known to depend
on the aspect angle as well as the position of the target \cite{Stakk2009}.
In some applications \cite{McIntyre1998}, it may be of interest 
to do Kalman filtering in polar coordinates or other coordinates
that induce a state dependence in measurement errors. 
Modeling of process error covariance may also be state dependent --- for example, 
Bar-Shalom \cite{YAA} suggests that the right choice of process noise level
for flight tracking models depends on the turn rate range expected. 
Therefore if we are estimating turn rate as part of the state, 
the process noise level can be modeled as a function of (a portion of) the state. 

These ideas motivate extensions of the standard Kalman smoothing formulation
to situations where process and measurement variances have known functional 
dependence on the state. Several such extensions have already been considered.
In \cite{Stakk2009}, the Unscented Transform is used to 
fit models with state-dependent matrices acting on 
observation noise. 
Linear systems with additive observation noise 
where measurement error variance is a known function of the state 
are studied in \cite{Zehnwirth1988}. 
Linear systems with 
control inputs transformed by state-dependent matrices are 
considered in \cite{Dutka2006}. Finally, adaptive
system identification, as presented in~\cite{chui2008kalman}, 
also falls into this class.

In this paper, we formulate the state-dependent 
covariance problem as a statistical estimation 
problem, and develop algorithms for obtaining
the {\it maximum a posteriori} (MAP) estimate. 
The ideas presented here extend those developed in \cite{Bell1995} 
for diagonal covariance matrices in kinetic tracer studies.
In the theoretical development, we allow the process 
and measurement functions to be nonlinear, 
and we allow the functional dependence of
covariance on the state to be nonlinear as well. 

The paper proceeds as follows. 
In Section \ref{StateDependence}, 
we review the statistical origins of the Kalman smoother, 
casting it as a structured nonlinear regression problem. 
We show that consideration of state-dependent variance 
in such a regression brings to the forefront terms that are usually ignored,
and develop an extended MAP objective to optimize. 
The proposed formulation can be used
for general nonlinear regression where
variance depends in a known functional way on the parameters. 
In Section \ref{ConvexComposite} we build 
a new algorithm for solving the 
resulting optimization problems, 
exploiting their {\it convex composite structure}. 
The key step is a special convex subproblem, which we 
solve in Section~\ref{sec:Newton}.
In Section~\ref{EGKF}, we show the necessary details 
required to implement this method for time series analysis, 
so as to preserve the computational complexity of the 
classic smoothing algorithms. In Section \ref{EGnumerics},
we provide a numerical experiment using simulated
data that demonstrates the performance of the new 
smoother and the potential modeling capabilities of the approach.  
{We end with conclusions. }


\section{Kalman Smoothing with State-Dependent Uncertainty}
\label{StateDependence}

The dynamic structure of the Kalman smoothing problem is specified as follows: 

\begin{equation}
\label{IntroGaussModel}
\begin{array}{rcll}
	\B{x_1}&= & g_1(x_0)+\B{w_1},
	\\
	\B{x_k} & = & g_k (\B{x_{k-1}})  + \B{w_k}& k = 2 , \ldots , N,
	\\
	\B{z_k} & = & h_k (\B{x_k})      + \B{v_k}& k = 1 , \ldots , N\;,
\end{array}
\end{equation}
where $g_k, h_k$ are known (nonlinear) process and measurement functions,
and $\B{w_k}\in\mB{R}^n$, $\B{v_k}\in\mB{R}^{m(k)}$ are mutually independent  Gaussian random variables
with positive definite covariance matrices $Q_k$ and $R_k$,  
$\B{x_k} \in \mB{R}^n$ are the unknown states, and $\B{z_k} \in \mB{R}^{m(k)}$ are the observed measurements. 

Considering model \eqref{IntroGaussModel} and using Bayes' theorem, the conditional likelihood 
of the entire state sequence $\{x_k\}$ given the measurement sequence $\{z_k\}$ is given by  
\begin{equation}\label{Bayes}
\B{p}\left(\{x_k\} \big| \{z_k\}\right) \propto \B{p}\left(\{z_k\}\big|\{x_k\}\right)\B{p}\left(\{x_k\}\right)\;, 
\end{equation}
which in turn can be written in terms of the likelihood of state increments $ \B{p}(w_k)$ and measurement residuals $\B{p}(v_k)$:
\begin{equation}
\label{MAP}
\begin{aligned}
&\B{p}\left(\{z_k\}\big|\{x_k\}\right)\B{p}\left(\{x_k\}\right) =
\kappa\prod_{k=1}^N  \B{p}(v_k)\B{p}(w_k) \\
& = 
\kappa\prod_{k=1}^N \exp\Big( -\frac{1}{2}(z_k - h_k(x_k))^\top R_k^{-1}(z_k - h_k(x_k)) \\
& \quad\quad-\frac{1}{2}(x_k - g_k( x_{k-1}))^\top Q_k^{-1}(x_k - g_k(x_{k-1}))\Big)\;,
\end{aligned}
\end{equation}
where we define $g_1(x_0) = x_0$. The constant of proportionality $\kappa$ is usually ignored, since in classic 
models, the variance terms $Q_k$ and $R_k$ are fixed. It is given by 
\begin{equation}
\label{VarConst}
\kappa = \prod_{k=1}^n \frac{1}{(2\pi)^{n/2} \det(Q_k)}\frac{1}{(2\pi)^{m(k)/2} \det(R_k)}.
\end{equation}

Our main contribution here is to remove the assumption that $Q_k$ and $R_k$ 
are fixed and known, and instead model 
these covariance matrices as known $\mathcal{C}^2$ functions of the state. 
In this setting, $\kappa$ in~\eqref{VarConst} is no longer a constant, and must be accounted for. 
To design our approach, we assume that 
we are given the inverse Cholesky factors $Q_k^{-1/2}(x_k)$ and $R_k^{-1/2}(x_k)$ 
as functions of the state.  For simple (e.g. diagonal variance) models, 
there is no loss of generality here; one can easily transform between different representations. 
When $Q_k$ and $R_k$ are full, however, the assumption that inverse Cholesky factors are available  
is essential for our approach. In addition to considerations of computational efficiency, 
the main motivation behind the assumption is the selection of an appropriate convex-composite model; 
this is explained in detail in the next section. 

In order to develop a simpler notation for estimating the entire 
state  sequence, we define
functions $g:\mathbb{R}^{nN}\rightarrow\mathbb{R}^{nN}$ and 
$h:\mathbb{R}^{nN} \rightarrow \mathbb{R}^{M}$, with $M = \sum_k m_k$,
from components $g_k$ and $h_k$ as follows:
\begin{equation}
\label{ghDef}
g(x) = \begin{bmatrix}x_1 \\ x_2 - g_2(x_1) \\ \vdots \\ x_N - g_N(x_{N-1}) \end{bmatrix}\;,
\quad
h(x) = \begin{bmatrix} h_1(x_1) \\ h_2(x_2) \\ \vdots \\ h_N(x_N)\end{bmatrix}\;.
\end{equation}

Given a sequence of column vectors $\{ u_k \}$
and matrices $ \{ T_k \}$ we use the following notation:
\begin{equation}\label{defs}
\begin{aligned}
&\R{vec} ( \{ u_k \} )
=
\begin{bmatrix}
u_1 \\ u_2  \\ \vdots \\ u_N
\end{bmatrix}
\; , \;
\R{diag} ( \{ T_k \} )
=
\begin{bmatrix}
T_1    & 0      & \cdots & 0 \\
0      & T_2    & \ddots & \vdots \\
\vdots & \ddots & \ddots & 0 \\
0      & \cdots & 0      & T_N
\end{bmatrix} , \\
&\begin{aligned}
R       & =  \R{diag} ( \{ R_k \} )
\\
Q       & =  \R{diag} ( \{ Q_k \} )
\\
x       & = \R{vec} ( \{ x_k \} )
\end{aligned}
\quad\quad\quad
\begin{aligned}
w      &  = \R{vec} (\{g_0, 0, \dots, 0\})
\\
z      & = \R{vec} (\{z_1,  \dots, z_N\}) 
\\
g_0 & = g_1(x_0).
\end{aligned} 
\end{aligned}
\end{equation}

With this notation, and under Gaussian assumptions,  the extended MAP object for the Kalman smoother, 
which incorporates state-dependent variance terms, is given by  
\begin{equation}
\label{fullNLLS}
\begin{aligned}
\frac{1}{2}\|Q^{-1/2}(x)(g(x) -w)\|_2^2 + \frac{1}{2}\|R^{-1/2}(x)(h(x) - z)\|_2^2\\
-\log\det\left(Q^{-1/2}(x)\right) - \log\det\left(R^{-1/2}(x)\right) \;.
\end{aligned}
\end{equation}
With~\eqref{fullNLLS} in front of us, we see why the log determinant terms play an important role.
Without these terms, an optimization approach to minimize a  weighted sum of squares will aim to drive 
$Q^{-1/2}(x)$ and $R^{-1/2}(x)$ to $0$ if at all possible. The function $-\log(\cdot)$ acts as a barrier to prevent 
this from happening.


\section{Convex Composite Formulation and Algorithm}
\label{ConvexComposite}

We would like to apply the generalized Gauss-Newton methodology for minimizing
convex composite functions \cite{Burke1985} to the objective~\eqref{fullNLLS}.
The first step in this process is to write this objective in convex composite form, that is, in the
form $f=\rho\circ F$, where $\rho$ is convex and $F$ is smooth. The choice of the functions
$\rho$ and $F$ depend on how we wish to model the representation of the problem.
The most straightforward way to rewrite \eqref{fullNLLS} is the more general form
\[
J(x):=\half c(x)^TW(x)^{-1}c(x)+\half\log\det(W(x)),
\]
where $\map{c}{\mB{R}^{nN}}{\mB{R}^{M+nN}}$ and $\map{W}{\mB{R}^{nN}}{\cS^{M+nN}_{++}}$
are smooth maps given by 
\begin{gather}
\label{c}
c(x) = \begin{bmatrix}x_1 - g_0\\h_1(x_1) -z_1\\ x_2 - g_2(x_1) \\ h_2(x_2) - z_2 \\ \vdots \\ x_N - g_N(x_{N-1})  \\ h_N(x_N) - z_N \end{bmatrix},\\
\label{W}
W(x) = \begin{bmatrix} Q_1(x_1) & 0 &&&\\ 0 & R_1(x_1) &&&\\ 
\vdots &\ddots &\ddots &&\vdots\\
&&& Q_N(x_N) &0 \\
&&&0&R_N(x_N)
\end{bmatrix}
\end{gather}
where  $\cS^{M+nN}_{++}$ is the cone of real symmetric $({M+nN})\times({M+nN})$
positive definite matrices.

Then $J=\hat\rho\circ \hat F$ with
\[
\begin{aligned}
\hat \rho(c,W)& :=\half c^TW^{-1}c+\half\log\det(W)\\
\hat F(x)&=(c(x),W(x)).
\end{aligned}
\]
Although the function $\hat F$ in this formulation can be assumed smooth, the function
$\hat\rho$ is not convex.  Indeed, $\hat\rho$ is the difference of two convex functions. 
When viewed as a function of $(c, W^{-1})$, it is still not jointly convex in these arguments.


Here, we propose an approach that applies in many practical
settings and yields an efficient solution procedure. 
However, a price is paid in a more complex model for the covariance matrices.
Specifically, 
we assume that 
the Cholesky factors for $Q_k^{-1}(x_k)$ and $R_k^{-1}(x_k)$ are given to us
as explicit functions of the state. 
We denote these factors by $Q_k^{-1/2}(x_k)$ and $R_k^{-1/2}(x_k)$,
respectively. In some settings, the matrices $Q_k(x_k)$ and $R_k(x_k)$ are modeled as
diagonal matrices, in which case the inverse Cholesky factors are easily computed diagonal matrices.
We provide an example of this type in the final section. Under this modeling assumption, 
the objective~\eqref{fullNLLS} can be abstracted to the more general form 
\begin{equation}
\label{extendedObjective}
K(x) 
=  
\frac{1}{2}c(x)^\R{T} V(x)^{T}V(x) c(x) 
- 
\log \circ \det [ V(x) ]\; ,
\end{equation}
where $\map{c}{\mB{R}^{nN}}{\mB{R}^{M+nN}}$ is exactly as in~\eqref{c} and 
$\map{V}{\mB{R}^{nN}}{\cL^{M+nN}}$ 
is given by 
\begin{equation}
\label{V}
V(x) = \begin{bmatrix} Q_1^{-1/2} & 0 &&&\\ 0 & R_1^{-1/2} &&&\\ 
\vdots &\ddots &\ddots &&\vdots\\
&&& Q_N^{-1/2} &0 \\
&&&0&R_N^{-1/2}
\end{bmatrix},
\end{equation}
where all blocks of V are functions of $x$, 
and $\cL^{M+nN}$ is the subalgebra of $({M+nN})\times({M+nN})$ real lower triangular matrices.
Throughout, we assume that
both $c$ and $V$ are twice continuously differentiable and that
$\mathrm{dom}(K):=\set{x}{K(x)<+\infty}=\set{x}{V(x)\in \cL^{M+nN}_{++}}\ne \emptyset$,
where $\cL^{M+nN}_{++}$ is the cone of $({M+nN})\times({M+nN})$ real lower triangular matrices
with strictly positive entries on the diagonal.
Now $K$ can be written in
convex composite form $K(x)= \rho \circ F$ with
%
%
\index{convex composite}
\begin{eqnarray}
\label{extended}
\rho(u, v) 
&=& \frac{1}{2}u^\R{T}u -\sum_{i}\log[ v_i] \\
\label{eq:define F}
\quad F(x)
&=&
\begin{bmatrix}
F_1(x)\\
F_2(x)\\
\end{bmatrix}
=
\begin{bmatrix}
V(x)c(x) \\
\R{vec}[\{V_{ii}(x)\}]
\end{bmatrix}\;.
\end{eqnarray}
Note 
 that $\mathrm{dom}(\rho)=\mB{R}^{M+nN}\times \mB{R}^{M+nN}_{++}$.

The direction finding subproblem in a Gauss-Newton method takes the form
\[
\min_d \rho(F(x)+F'(x)d) +\sfrac{\omega}{2} d^\R{T}d\ ,
\]
for some $\omega\ge 0$. The quadratic term $\frac{\omega}{2} d^\R{T}d$
is a regularization term that both guarantees the uniqueness of the solution 
and regulates its magnitude.
The convergence analysis of methods of this type rely heavily on 
the difference function 
\begin{equation}
\label{extendedDelta}
\Delta( x; d) = \rho \left( \; F(x) + F' (x) d \; \right)
 - K(x) \;.
\end{equation}
which is important for both convergence criteria
and the line search in the overall method. 
In particular, \cite[Lemma 2.3]{Burke1985}
\begin{equation}\label{eq:cvx comp dd}
K'(x:d)=\inf_{t>0} t^{-1}\Delta( x; td) \mbox{ for all }d\in\mB{R}^{nN},\; x\in\mathrm{dom}(K)
\end{equation}
since, whenever $F_2(x)>0$,  then, for all $d$, $F_2(x)+F'_2(x)(td)>0$ for all $t$ sufficiently small.

Linearizing 
the functions $F_i(x)$ in \eqref{extended} 
yields approximations
$\tilde F_i(x;d):=F_i(x)+F'_i(x)d$, which in turn gives the approximation 
\begin{equation}
\label{approxEO}
\tilde K(x;d) 
= 
\rho[\tilde F_1(x;d), \tilde F_2(x;d)]\; .
\end{equation}
This is the objective for the direction finding subproblem.
Here,
\begin{equation}
\label{approxF}
\begin{array}{lll}
\tilde F_1(x; d)
&=& 
\left(
V(x)\partial_x c(x) + (c(x)^\R{T}\otimes I_N)\partial_x V(x) \right) d\\
&&\quad + V(x)c(x)\\
\tilde F_2(x; d) 
&=&  
\R{vec}\left(\{V_{ii}(x) + \partial_x V_{ii}(x)d\}\right).
\end{array}
\end{equation}
Note that we must be sure that $\tilde F_2(x;d)$
is component-wise greater than zero.  For details 
of these derivations, see~\cite{AravkinThesis2010}.
The Gauss-Newton subproblem is now given by
\begin{equation}
\label{extendedGNSubproblem}
\begin{aligned}
&\bDel(x):=\min_{d\in \mB{R}^{nN}} \Delta( x; d)+\sfrac{\omega}{2} d^\R{T}d\\
&=\min_{d\in \mB{R}^{nN}}
	\frac{1}{2}\tilde F_1(x;d)^\R{T}\tilde F_1(x;d) +\sfrac{\omega}{2}d^\R{T}d
        - 
        \sum_{i} \log [\tilde F_2(x; d)]\;. 
        \end{aligned}
\end{equation}
Due to our assumptions on $c$ and $V$, these subproblems are always well defined,
are convex, and have a unique solution which must always exist. In addition, they provide an estimate
for the first-order optimality for $K$.
\begin{theorem}\cite[Theorem 3.6]{Burke1985}\label{thm:first-order opt}
Let $x\in\mathrm{dom}(K)$. Then the 
following three statements are equivalent:
\begin{enumerate}
\item[(i)] $\bDel(x)=0$,
\item[(ii)] $\bar d=0$ solves \eqref{extendedGNSubproblem}, and
\item[(ii)] $0\in\partial K(x)$, where $\partial K(x)=F'(x)^T\partial \rho(F(x))$ is the 
generalized subdifferential of $K$ at $x$ \cite[Definition 8.3]{RTRW}.
\end{enumerate}
In particular, these conditions imply that $x$ is a first-order stationary point for $K$.
\end{theorem}

If we ignore dependence on $x$,
the optimization problem in \eqref{extendedGNSubproblem} can be rewritten as 
\begin{equation}
\label{extendedGNSubproblemTwo}
\begin{aligned}
\displaystyle
\min_{d \in \mB{R}^{nN}} 
	 \quad&\frac{1}{2}d^\R{T}Cd + a^\R{T}d 
        - 
        \sum_{i} \log [s_i] \\
   \mbox{s.t.} 
   \quad & \quad s = \R{vec}\{V_{ii}\} + \partial_x \R{vec}\{V_{ii}\}d ,
\end{aligned}
\end{equation}
where
\begin{equation}
\label{extendedHessian}
\begin{aligned}
C &= \omega I+\widetilde{V}^\R{T} \widetilde{V}\\
a &= \widetilde{V}^\R{T}Vc\\
\widetilde{V}&=[V\partial_xc\! +\! (c^\R{T}\otimes I_N)\partial_x V] 
\end{aligned}
\end{equation}
Note that $a$ is the gradient of the quadratic portion of the extended objective with respect to 
the state sequence $x$. The quantity $(c^\R{T}\otimes I_N)\partial_x V$
that appears in \eqref{extendedHessian} can be rewritten as
\begin{equation}
\begin{array}{lll}
(c^\R{T}\otimes I_N)\partial_x V
&=&
\sum_{i=1}^{M + nN} c_i \partial_x V_{i\cdot}\ .
\end{array}
\end{equation}
The Lagrangian associated with the extended subproblem~\eqref{extendedGNSubproblemTwo} 
is given by  
\begin{equation}
\label{extendedLagrangian}
\begin{aligned}
L(d, s, \lambda) 
&=  
\frac{1}{2}d^\R{T}Cd + a^\R{T}d 
- 
\sum_{i} \log [s_i]\\
&+ \lambda^\R{T} \left(s - \R{vec}\{V_{ii}\} - \partial_x \R{vec}\{V_{ii}\}d\right),
\end{aligned}
\end{equation}
for $s>0$ and $\lambda >0$.
The corresponding optimality conditions state that a direction $d$ 
solves~\eqref{extendedGNSubproblemTwo}
if and only if there exist 
$s,\lambda\in\mB{R}^{M+nN}_{++}$ such that
\begin{equation}
\label{extendedOptimalityConditions}
\begin{array}{lll}
\nabla_{d}L &=& Cd + a - \partial_x \R{vec}\{V_{ii}\}^\R{T}\lambda = 0\\
\nabla_{s}L &=& -D(s)^{-1}\B{1} + \lambda = 0\\
\nabla_{\lambda} L &=& s - \R{vec}\{V_{ii}\} - \partial_x \R{vec}\{V_{ii}\}d = 0\; ,
\end{array}
\end{equation}
where $D(s):=\diag(s)$.

We refer to~\eqref{extendedGNSubproblemTwo} as the {\it extended subproblem}. 
In the next section, we show that this problem can be rapidly solved. 
This motivates the {\it Extended Gauss-Newton} method for~\eqref{extendedObjective}.

\begin{algor}
\label{GaussNewtonAlgorithm}
{\it Generalized Gauss-Newton Algorithm.}

The inputs to this algorithm are
\begin{itemize}
\item \( x^0 \in \mathrm{dom}(K):=\{x:K(x)<\infty\}\subset\mB{R}^{Nn} \): initial estimate of state sequence
\item \(\varepsilon \ge 0 \): overall termination criterion
\item \(\omega > 0 \): regularization parameter
\item \( \beta \in (0, 1) \): step size selection parameter
\item \( \gamma \in (0,1) \): line search step size factor
\end{itemize}

The steps are as follows:

\begin{enumerate}
\item
Set the iteration counter \( \nu = 0 \).
\item (Generalized Gauss-Newton Step)
\label{GaussNewtonStep}

Find descent direction \(d^\nu \)
solving \eqref{extendedGNSubproblemTwo} and set 
$\Delta_\nu := \bDel(x^\nu)=\Delta( x^\nu ; d^\nu)$.
{\it Terminate} if \( \Delta_\nu \ge - \varepsilon \).
\item (Line Search) Set
\[
\begin{array}{lll}
t_\nu &=& \max  \gamma^i \\
     &\text{s.t.}&i \in \{ 0, 1, 2, \cdots \} \; \mbox{ and }
    \\
    &\text{s.t.}&\rho \left( F( x^\nu + \gamma^i d^\nu ) \right)
        \le \rho \left( F( x^\nu ) \right)
                     + \beta \gamma^i \Delta_\nu.
\end{array}
\]
\item (Iterate)
Set \(x^{\nu+1} = x^\nu + t_\nu d^\nu \) and return to Step~\ref{GaussNewtonStep}.
\end{enumerate}
\end{algor}

\begin{remark}
Note that the line search is well defined whenever $\Delta_\nu\ne 0$. 
Indeed, 
since whenever $\mathrm{diag}(V(x))>0$, then,
for all $d$, $F_2(x;t d)>0$ for all $t$ sufficiently small.
Consequently, since $x^0\in\mathrm{dom}(K)$, we have $\{x^\nu\}\subset\mathrm{dom}(K)$.
In addition, by \eqref{eq:cvx comp dd}, $K'(x^\nu;d^\nu)\le\Delta_\nu <\bDel(x^\nu)$,
so that $\gamma^{-i}\Delta(x^\nu;\gamma^{i}d^\nu)<\beta \Delta_\nu$ for all $i$ sufficiently large.
\end{remark}


\begin{theorem}
\label{thm:cvg}[Convergence]
Let $\{x^\nu\}$ be generated by Algorithm \ref{GaussNewtonAlgorithm} with $\epsilon = 0$.
Then either the algorithm terminates finitely at a first-order stationary point for $K$
or the sequence $\{x^\nu\}$ is infinite and every cluster point of the sequence is a 
first-order stationary point for $K$.
\end{theorem}
The proof is given in the appendix. 

In the next section, we show how to solve the subproblem~\eqref{extendedGNSubproblemTwo} 
for a general state-dependent covariance
regression problem.



\section{Solving the Extended Subproblem}
\label{sec:Newton}

%
To solve the direction finding subproblem,
we apply a damped Newton method directly to the optimality conditions
\eqref{extendedOptimalityConditions}. 
We present the high-level method here, with details concerning Kalman smoothing
given in the next section. 

Let $E(s, \lambda, d)$ denote the KKT system given in~\eqref{extendedOptimalityConditions}, 
rearranged in a particular order:
\begin{equation}
\label{compE}
E(s, \lambda, d) 
=
\begin{bmatrix} 
s - \R{vec}\{V_{ii}\} - \partial_x \R{vec}\{V_{ii}\}d\\
D(s)D(\lambda)\B{1} - \B{1} \\
Cd + a - \partial_x \R{vec}\{V_{ii}\}^\R{T}\lambda
\end{bmatrix}\\ 
\end{equation}

Our goal is to find $(\overline s, \overline\lambda, \overline d)$ 
for which $E(\overline s, \overline\lambda, \overline d)=0$.
We use
damped Newton's method on E which requires solving the Newton equation 
\begin{equation}
\label{eq:Newton}
\nabla E(s, \lambda, d) 
\begin{bmatrix}
\Delta s\\
\Delta \lambda\\
\Delta d
\end{bmatrix}
= 
-E(s, \lambda, d),
\end{equation}
where $\nabla E(s, \lambda, d) $ is given by 
\begin{equation}
\label{compEgrad}
\begin{bmatrix}
I & 0 & -\partial_x \R{vec}\{V_{ii}\} \\
D(\lambda) & D(s) & 0 \\
0 & -\partial_x \R{vec}\{V_{ii}\}^\R{T} &C
\end{bmatrix}.
\end{equation}
Define 
\[
\Sc{V} = \partial_x \R{vec}\{V_{ii}\}.
\]
Then, using row operations 
\[
\begin{aligned}
R_2 &= R_2 - D(\lambda)R_1\\
R_3 &= R_3 + \Sc{V}^\R{T}D(s)^{-1}R_2,
\end{aligned}
\] 
we obtain the modified system 
\begin{equation}
\begin{array}{lll}
&&\begin{bmatrix}
I & 0 & -\Sc{V}\\
0 & D(s) & D(\lambda)\Sc{V}\\
0 & 0 & \Phi
\end{bmatrix}
\begin{bmatrix}
\Delta s\\
\Delta \lambda\\
\Delta d
\end{bmatrix}
=
\begin{bmatrix}
\alpha\\
\beta\\
\gamma
\end{bmatrix}\;, 
\end{array}
\end{equation}
where
\begin{equation}
\label{eq:Phi}
\Phi = C + \Sc{V}^\R{T}D(s)^{-1}D(\lambda)\Sc{V}
\end{equation}
and 
\[
\begin{array}{lll}
 \alpha &=& -s + \R{vec}\{V_{ii}\} + \Sc{V}d\\
 \beta &= &\B{1} - D(\lambda)\left(\R{vec}\{V_{ii}\} + \Sc{V}d\right)\\
\gamma &=&  
\Sc{V}^\R{T}
\Big(
\lambda 
+ 
D(s)^{-1}
\left(
\B{1} - D(\lambda)(\R{vec}\{V_{ii}\} + \Sc{V}d)
\right)
\Big) \\
&&-Cd - a 
\;.
\end{array}
\]
By \eqref{extendedHessian}, the matrix $C$ is always positive definite and so $\Phi$ is
always positive definite and hence invertible.
This allows us to recover the Newton direction: 
\begin{equation}
\label{extendedDirection}
\begin{array}{lll}
\Delta d 
&=& 
\Phi^{-1}\gamma\\
\Delta \lambda 
&=&
D(s)^{-1}\Big(\B{1} - D(\lambda)(\R{vec}\{V_{ii}\} + \Sc{V}(d + \Delta d))\Big)\\
\Delta s 
&=&
-s + \R{vec}\{V_{ii}\} + \Sc{V}(d + \Delta d)\ .
\end{array}
\end{equation}
Note that any damping scheme requires that 
$s > 0$ for the objective to be finite, and 
hence, in addition, we require that $\lambda > 0$ since we need $D(s)D(\lambda)\B{1}=\B{1}$ . 

\begin{figure*}
\begin{center}
{\includegraphics[scale=0.7]{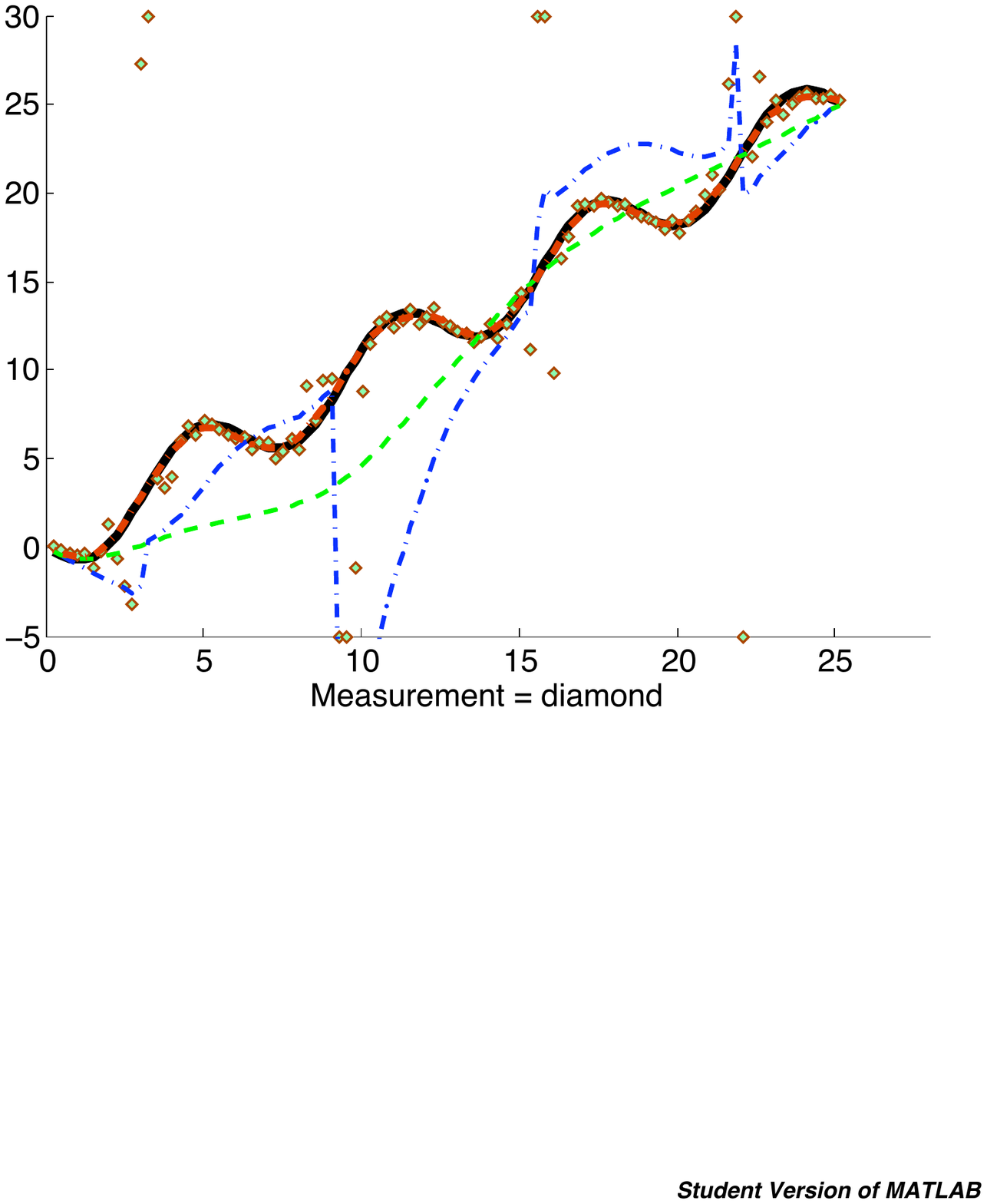}}
\end{center}
\caption{\label{ExtendedGN}
True state $x_1$ (black curve), Extended Smoother estimate (thick red dash-dot), Kalman filter estimate (blue dash-dot) and Kalman Smoother estimate (green dashed curve). Measurements are displayed as diamonds, and those outside the axis range are displayed on the figure boundary. 
}
\end{figure*}

\section{Structure of the Extended Kalman Smoothing Objective}
\label{EGKF}

We now specify the method in the previous section to the Kalman smoothing problem,
and demonstrate that the computational efficiency of the Kalman smoother can be preserved. 

The functions $c(x)$ and $V(x)$ are given by~\eqref{c} and~\eqref{V}.

With these definitions, objective $K(x)$
in~(\ref{extendedObjective}) is exactly~\eqref{fullNLLS}, 
and can be written explicitly as follows: 
\begin{equation}
\label{explicitEKS}
\begin{array}{lll}
&&\frac{1}{2}\Big(c^\R{T}(x)V(x)^\R{T}V(x)c(x)\Big)  - \log\circ\det[ V(x)]\\
& =  &
\displaystyle
\frac{1}{2}\sum^N_{k=1}  
\|[z_k - h_k(x_k)]\|^2_{R_k^{-T/2}(x_k)R_k^{-1/2}(x_k)} \\
&&+
\displaystyle
\frac{1}{2}\sum^N_{k=0} 
\|x_k - g_k(x_{k-1})\|^2_{Q_k^{-T/2}(x_k)Q_k^{-1/2}(x_k)}\\
&&  -\log\det( R^{-1/2}_k(x_k))
- \log\det( Q^{-1/2}_k(x_{k})) ,
\end{array}
\end{equation}
where, for any symmetric positive definite matrix $Q$, $\|u\|^2_Q:=u^TQu$.

We now derive the explicit forms for $C$ and $a$ in (\ref{extendedHessian}) 
for the Gauss-Newton subproblem (\ref{extendedGNSubproblemTwo}). Recall
that $C$ and $a$ are given by 
\[
\begin{aligned}
 C& =\! \omega I\!+\!
[V\partial_xc \!+ \!(c^\R{T}\otimes I_N)\partial_x V]^\R{T} 
[V\partial_xc \!+\! (c^\R{T}\otimes I_N)\partial_x V], \\
a &=\! c^\R{T}V^\R{T} V\partial_x c 
+
c^\R{T}V^\R{T} (c^\R{T}\otimes I_N)\partial_x V.
\end{aligned}
\]
where
\[
\begin{array}{lll}
(c^\R{T}\otimes I_N)\partial_x V
&=&
\sum_{i=1}^{M + nN} c_i \partial_x V_{i\cdot},  
\end{array}
\]
and 
\begin{equation}
\label{eq:full-c}
\partial_x c(x) = 
\begin{bmatrix}
    \R{I}  & 0      &&          &
    \\
            H_1 & 0 && \dots & 0
        \\
    -G_2   & \R{I}  && \ddots   &
    \\
            0 & H_2 & 0 && \dots
        \\
        & -G_3 &  \ddots  && 0
    \\     
        &  &  \ddots  &H_{N-1}& 0
    \\
        &        &&   -G_N  & \R{I}
        \\
        & & &0 & H_N 
\end{bmatrix}\;,
\end{equation}
with $G_k = \partial_{x_k} g_{k+1}(x_k)$, $H_k = \partial_{x_k} h_k(x_k)$, 
and the dependence on $x$ has been suppressed to decrease the notational burden. 
Note that the matrix $G$ is invertible, and so $\partial_x c(x)$ is injective, that is,
$\mathrm{Null}(\partial_x c(x))=\{0\}$. In addition, since, we require 
$\R{vec}\{V_{ii}\}>0$ at every iteration, the matrix $V^TV$ is always positive definite.
Consequently, the matrix $ \partial_x c^TV^TV\partial_x c$ is always positive definite.

Define $\tilde w(x)$ and $\tilde v(x)$ in~\eqref{defs} by 
\begin{eqnarray}
\label{wk} 
\tilde w_k(x) &=& x_k - g_k(x_{k-1})\\
\label{vk}
\tilde v_k(x) &=& z_k - h_k(x_k)
\end{eqnarray}
In the expressions below, we will use notation $\tilde w_{k,i}$ to mean
the $i$th component of $\tilde w_k$. 

The matrix $(c^\R{T}\otimes I_N)\partial_x V$ has the following block structure:
\begin{equation}
\label{Vderivative}
[(c^\R{T}\otimes I_N)\partial_x V] 
=
\small
\begin{bmatrix} 
\tilde Q_1 & 0 & 0 & 0 &0 & \cdots & 0 \\
\tilde R_1 & 0 & 0 & 0 &0 & \cdots & 0 \\
0 & \tilde Q_2 & 0 & 0 &0 & \cdots & 0 \\
0 & \tilde R_2 & 0 & 0 &0 & \cdots & 0 \\
0 & 0 & 0 & \ddots & \ddots & \cdots & 0 \\
0 & 0 & 0 & \ddots & \ddots & Q_{N-1} & 0\\
0 & 0 & 0 & \ddots & \ddots & R_{N-1} & 0\\
0 & 0 & 0 & \cdots & \cdots & 0 &\tilde Q_N \\
0 & 0 & 0 & \cdots & \cdots & 0 &\tilde R_N \\
\end{bmatrix}
\end{equation}
where 
\begin{eqnarray}
\label{Qderivative}
\tilde Q_i 
&=&
\sum_{j=1}^n \tilde w_{i,j} \partial_{x_i}[Q_i^{-1/2}]_{j\cdot}\\
\label{Rderivative}
\tilde R_i 
&=& 
\sum_{j=1}^{m(i)} \tilde v_{i,j}\partial_{x_i}[R_i^{-1/2}]_{j\cdot}\;.
\end{eqnarray}
Then we can write down the gradient $a$: 
\[
\begin{array}{lll}
a 
&=&
\begin{bmatrix}
a_1^\R{T} & a_2^\R{T} & \cdots & a_N^\R{T}
\end{bmatrix}^\R{T}\;,
\end{array}
\]
where
\begin{equation}
\label{extendedGradient}
\begin{array}{lll}
a_j
&=&
-
\tilde v_j^\R{T}R^{-1}_j(x_{j})\partial_{x(j)}h_j(x_j) 
+  
\tilde w_j^\R{T}Q^{-1}_j(x_{j-1})\\
&&+
\tilde w_{j+1}^\R{T}
Q^{-1}_{j+1}(x_{j})
\partial_{x(j)}g_{j+1}(x_{j}) \\
&&+
\tilde w_j^\R{T}Q_{j}^{-T/2}\tilde Q_j 
+
\tilde v_j^\R{T}R_j^{-T/2}\tilde R_j\;.
\end{array}
\end{equation}

Using~\eqref{eq:full-c} and~\eqref{Vderivative}, we can form 
\[
\Psi := V\partial_x c + (c^\R{T}\otimes I_N)\partial_x V,
\]
%
%
and obtain a closed form solution for $C = \Psi^T\Psi$:
\begin{equation}
\label{extendedC}
\begin{aligned}
C
&=
\begin{bmatrix} C_1 & A^\R{T}_2 & 0 &  \\
A_2 & C_2 & A^\R{T}_3 & 0 \\
0 & \ddots & \ddots & \ddots\\
 & 0 & A_N &C_N
\end{bmatrix}\\
C_k 
&\!=\!\omega I\!+\!
[Q_k^{-1/2}\!\! +\! \tilde Q_k]^\R{T}[Q_k^{-1/2}\!\! +\! \tilde Q_k]
\!+\!
\nabla g_{k+1}^\R{T} Q_{k+1}^{-1}\!\nabla g_{k+1}\\
&\ \ +
[-R_k^{-1/2}\nabla h_k + \tilde R_k]^\R{T}
[-R_k^{-1/2}\nabla h_k + \tilde R_k]\\
A_k
&=
-(Q_{k}^{-1/2}+\tilde Q_k)^\R{T}Q_k^{-1/2}\nabla g_k\;.
\end{aligned}
\end{equation}

\begin{remark}
\label{rem:implementable}
The matrix (\ref{extendedC}) 
is block tridiagonal, and so it can be inverted with effort $O(n^3N)$ using
any of the algorithms in~\cite{Bell2000,AravkinBellBurkePillonetto2013}. 
The proof sketch is given in the appendix. 
\end{remark}

Recall the direction finding equation (\ref{eq:Newton})
in the Extended GN algorithm:

\begin{equation*}
\nabla E(s, \lambda, d) 
\begin{bmatrix}
\Delta s\\
\Delta \lambda\\
\Delta d
\end{bmatrix}
= 
-E(s, \lambda, d). 
\end{equation*}
The solution to this system is given by (\ref{extendedDirection}), 
with $\Phi$ as in~\eqref{eq:Phi}, $C$ as in~\eqref{extendedC}, and 
$[\Sc{V}^\R{T}D(s)^{-1}D(\lambda)\Sc{V}]_k$ 
given by 
\[
\begin{array}{lll}
\partial_{x(k)}\R{diag}\{Q_{k}^{-1/2}\}^\R{T}
D(s_{Q_k})^{-1}D(\lambda_{Q_k})
\partial_{x(k)}\R{diag}\{Q_{k}^{-1/2}\}\\
+
\partial_{x(k)}\R{diag}\{R_{k}^{-1/2}\}^\R{T}
D(s_{R_k})^{-1}D(\lambda_{R_k})
\partial_{x(k)}\R{diag}\{R_{k}^{-1/2}\}\;.
\end{array}
\]

\section{Numerical Results}
\label{EGnumerics}

\index{numerical results!extended Kalman smoother}
In this section, we present some numerical experiments to show 
the advantages and modeling possibilities of the new Kalman smoother.
The simulation model we consider is similar to the one presented
in~\cite{Bell2008}.
The `ground truth' time series
for this simulated example is given by 
\[
x(t) = \begin{bmatrix}1 -2\cos(t) \\ t - 2\sin(t)\end{bmatrix} \; .
\]
The time between measurements is a constant
denoted by $\Delta t$.
The models for the mean of $x_k$ given $x_{k-1}$ and for process covariance $Q_k$~\cite{Jaz,Oks} are 
\[
g_k ( x_{k-1} ) =
\begin{bmatrix}
	1        & 0
	\\
	\Delta t & 1
\end{bmatrix} x_{k-1}
\;, \quad
Q_k =
\begin{bmatrix}
	\Delta t   & \Delta t^2 / 2
	\\
	\Delta t^2/2 & \Delta t^3 / 3
\end{bmatrix}.
\]
The measurement model for the mean of $z_k$ given $x_k$ is
\(
h_k( x_k )  = x_{2,k} \; ,
\)
where \( x_{2,k} \) denotes the second component of \( x_k \).

The main innovation of the example is in the measurement variance model. 
The smoother takes inverse Cholesky factors as input, 
and these are assumed to be $3 - x_{1, k}$. 
Then the variance model is given by 
$R_k(x_k) = (3 - x_{1,k})^{-2}$.  The measurements 
were generated using the measurement model,
from two full periods of the time series $x(t)$, with \( N = 100 \)
discrete time points equally spaced over the interval $[0, 4 \pi]$, 
and  with noise sampled from $N(0, R_k(x_k))$.
Since the true state $x_1$ varies in the interval $[-1,3]$, 
the variance for the observations goes to infinity 
when $t$ is a multiple of $\pi$. 

This simulation illustrates a situation 
where the measurements are very reliable for some state values, 
but completely unpredictable for others. This phenomenon
may occur for example if sensors report garbage values when the attitude of a vehicle
is in a particular configuration. The measurement model
presented here can be easily adapted by the user to take their 
beliefs about the system into account. 
The main point is that as long as the inverse Cholesky factors for the 
variance can be coded
as a smooth function of the state, smoothed 
estimates for state values can be obtained taking into account 
this bad behavior of the measurements. 
 
The result of the simulation is shown in Figure \ref{ExtendedGN}.
The extended Kalman smoother (thick red dash-dot) is able to recover 
the ground truth (shown in black) with no appreciable difference. 
The Kalman filter (thin blue dash-dot) 
is strongly affected by the outlying measurements, 
as expected. The Kalman smoother (green dashed) is able
to smooth the measurements, but cannot pick up the oscillations
of the ground truth, which are small in magnitude compared 
to the size of the errors.

This last point is the most important --- it is not just
the magnitude of the outliers that makes the Kalman smoother
fail, although it can be seen to be rather far off the 
ground truth. The biggest challenge of the situation presented is 
knowing which measurements to trust, since this information
depends on the state being estimated. 

\section{Conclusions}
In this paper, we presented extended formulations for modeling 
dynamic systems in cases where the covariance matrices are known
functions of the state. The formulation includes variance-control terms
arising from statistical modeling assumptions. These terms give rise  
to an extended convex-composite structure, and we propose a 
new method, the extended Gauss-Newton, which repeatedly solves 
extended convex subproblems by exploiting their KKT optimality conditions. 
When applied to dynamic inference problems, the proposed approach 
preserves the complexity of the classic Kalman smoother.  

\bibliographystyle{plain}
\bibliography{filter}

\section{Appendix}
\subsection{Proof sketch of Remark~\ref{rem:implementable}}
Recall the matrix inversion lemma:
\begin{lemma}[Matrix Inversion Lemma]
\label{lem:woodbury} Assume matrices $M_1 \in \mathbb{R}^{m_1\times m_1}$
and $M_2 \in \mathbb{R}^{m_2 \times m_2}$ are symmetric positive definite. Then 
for any matrix $U \in \mathbb{R}^{m_1 \times m_2}$, the following matrix is 
also symmetric positive definite:
\begin{equation}
\label{eq:inverse}
S = M_1^{-1} - M_1^{-1} U(M_2 + U^T M_1^{-1} U)^{-1}U^TM_1^{-1}.
\end{equation}
\end{lemma}
From this, we get an immediate and useful corollary:
\begin{corollary}
\label{cor:implementable}
For a positive definite matrix $M_2$, and any matrix $U$,
we have  
\[
\|U^T(M_2+UU^T)^{-1}U\|_2 < 1,
\]
where $\|\cdot\|_2$ denotes the spectral norm. 
\end{corollary}
Simply take $M_1 = I$ in lemma~\ref{lem:woodbury}. 
The conclusion of the lemma gives the result. 

Corollary~\ref{cor:implementable} can be applied to show that 
the algorithms in~\cite{AravkinBellBurkePillonetto2013} 
yield invertible blocks at every iteration.
Application of these algorithms 
to $C$ in~\eqref{extendedC} requires inverting matrices of form 
\[
H^+ + V^T(I - U^T(H + UU^T)^{-1}U)V,
\]
where $H$ and $H^+$ are always positive definite. 
The second term in the sum is clearly seen to be positive semidefinite 
by Corollary~\ref{cor:implementable}.
To be specific, at the first iteration, 
\[
\begin{aligned}
V &= Q_2^{-1/2} + \tilde Q_2, \\
H &= \omega I + (Q_1^{-1/2} + \tilde Q_1)^T(Q_1^{-1/2} + \tilde Q_1) \\
&+ (-R_1^{-1/2}\nabla h_1 + \tilde R_1)^\R{T}(-R_1^{-1/2}\nabla h_1 + \tilde R_1),\\
U &= \nabla g_2^T Q_2^{-1/2},\\
H^+ & = \omega I + (-R_2^{-1/2}\nabla h_2 + \tilde R_2)^\R{T}(-R_2^{-1/2}\nabla h_2 + \tilde R_2)
\end{aligned}
\]
The full argument can be made by induction, but we do not include it here. 

\subsection{Proof of Theorem~\ref{thm:cvg}}
The algorithm can only terminate if 
\[
0=\Delta_\nu=\Delta(x^\nu;d^\nu)\le \Delta(x^\nu;d^\nu)+\sfrac{\omega}{2}\|d^\nu\|^2=\bDel(x^\nu)\le 0,
\]
i.e., $\bDel(x^\nu)=0$, or equivalently, $x^\nu$ is a first-order stationary point for $K$
by Theorem \ref{thm:first-order opt}. 

Assume that the algorithm does not terminate finitely, and let $\hx$ be a cluster point
of the sequence of iterates $\{x^\nu\}$. Since this is a descent algorithm, it is necessarily
the case that $K(x^\nu)\downarrow K(\hx)$. 
Let $J\subset\mB{N}$ and $\hx\in\mB{R}^{nN}$ be such that 
$x^\nu\overset{J}{\rightarrow} \hx$, and suppose to the contrary that 
$\hx$ is not a first-order staionary point for $K$, i.e. $\bDel(\hx)< 0$.
We now use the optimality conditions \eqref{extendedOptimalityConditions} to show that
the subsequence of search directions $\{d^\nu\}_J$ is bounded. 
Let $(s^\nu,\lambda^\nu,d^\nu)$ denote the triple satisfying these conditions for each
$x^\nu$. Then multiplying the second condition in \eqref{extendedOptimalityConditions}
by $s^\nu$ 
and the third condition in \eqref{extendedOptimalityConditions} by $\lambda^\nu$
and combining, we find that
\[
M+nN=(\lambda^\nu)^\R{T}(\R{vec}\{V_{ii}(x^\nu)\} + \partial_x \R{vec}\{V_{ii}(x^\nu)\}d^\nu).
\]
By combining this with the first condition in \eqref{extendedOptimalityConditions}, we find that
\begin{equation*}
\begin{aligned}
M+nN&=(\lambda^\nu)^\R{T}\R{vec}\{V_{ii}(x^\nu)\}+(C(x^\nu)d^\nu+a(x^\nu))^\R{T}d^\nu\\
&\ge a(x^\nu)^\R{T}d^\nu+(d^\nu)^\R{T}C(x^\nu)d^\nu\\
&\ge\omega \Vert d^\nu\Vert^2_2 -\Vert a(x^\nu)\Vert_2\Vert d^\nu\Vert_2,
\end{aligned}
\end{equation*}
where the first inequality follows since $\lambda^\nu>0$ and $\diag(V(x^\nu))>0$,
and the second inequality follows from \eqref{extendedHessian}. 
Consequently, the subsequence $\{d^\nu\}_J$ is bounded due to the continuity
of 
\[
a(x)=\nabla\left(\half c(x)^\R{T}V(x)^\R{T}V(x)c(x)\right).
\]

With no loss in generality, we can now assume that there is a $\hat d$ such that
$d^\nu\overset{J}{\rightarrow}\hat d$. By continuity,
\[
\Delta_\nu=\Delta(x^\nu;d^\nu)\rightarrow \Delta(\hx;\hat d).
\] 
Moreover, for all 
$d\in\mB{R}^{M+nN}$,
\[
\Delta(x^\nu;d^\nu)+\sfrac{\omega}{2} (d^\nu)^\R{T}d^\nu
\le \Delta(x^\nu;d)+\sfrac{\omega}{2} d^\R{T}d\ .
\]
Taking the limit over $\nu\in J$ gives
\[
\Delta(\hx;\hat d)+\sfrac{\omega}{2} \hat d^\R{T}\hat d
\le \Delta(\hx;d)+\sfrac{\omega}{2} d^\R{T}d\ .
\]
Therefore, $\Delta(\hx;\hat d)+\sfrac{\omega}{2} \hat d^\R{T}\hat d=\bDel(\hx)$.

Recall our working assumption that
$\bDel(\hx)<0$.  Since we have just shown that $\Delta_\nu\rightarrow \bDel(\hx)$, 
we must therefore have $\xi:=\sup_{\nu\in J}\Delta_\nu<0$. 
Since
$K(x^{\nu +1})-K(x^\nu)\le \beta t_\nu\Delta_\nu$ with $K(x^\nu)$ convergent,
we must have $t_\nu\rightarrow 0$. Again, with no loss in generality,
$1> t_\nu\downarrow_J 0$. 
By continuity, there are $\delta>0$ and $\mu>0$ such that 
\[
\begin{aligned}
&\diag(V(x))\in \diag(V(\hx))+\mu\mB{B}\subset \mB{R}^{M+nN}_{++}\ \mbox{ and }\\
&\diag(V(x) +V'(x)d)\in \diag(V(\hx))+\mu\mB{B}
\end{aligned}
\]
for all $x\in\hx +\gamma^{-1}\delta \mB{B}$ and $d\in \gamma^{-1}\delta\mB{B}$.
Since $\{d^\nu\}_J$ is bounded and $t_\nu\downarrow_J0$, we can assume with
no loss in generality that $x^\nu\in \hx+\delta\mB{B}$ and 
$t_\nu d^\nu\in\delta\mB{B}$ with $1>t_\nu$ for all $\nu\in J$. 
Let 
\[
\kappa_1:=\sup\left\{\half\Vert V(x)c(x)\Vert_2^2\; :\; x\in \hx +\delta \mB{B}\right\}.
\]
Since
$\kappa_1\mB{B}\times \mu\mB{B}\subset \mathrm{intr}(\mathrm{dom}(\rho))$, 
$\rho$ is Lipschitz continuous on $\kappa_1\mB{B}\times \mu\mB{B}$ with Lipschitz
constant $\kappa_2>0$. Also, the function $F$ defined in \eqref{eq:define F} is such that
$F'$ is Lipschitz continuous on $\hx +\gamma^{-1}\delta \mB{B}$ with Lipschitz constant $\kappa_3>0$.

Due to the way the step sizes $t_\nu$ are chosen and the fact that 
$t_\nu d^\nu\in\delta\mB{B}$ and $1> t_\nu$ for all $\nu\in J$, we have
\[
\begin{aligned}
&\gamma^{-1}t_\nu\beta\bDel(x^\nu)< K(x^\nu+\gamma^{-1}t_\nu d^\nu)-K(x^\nu)\\
&\le \Delta(x^\nu;\gamma^{-1}t_\nu d^\nu)\\
&\quad+\kappa_2\Vert F(x^\nu+\gamma^{-1}t_\nu d^\nu)-F(x^\nu)-F'(x^\nu)(\gamma^{-1}t_\nu d^\nu)\Vert\\
&\le \gamma^{-1}t_\nu
\bDel(x^\nu)+\frac{\kappa_2\kappa_3}{2}(\gamma^{-1}t_\nu)^2\Vert{d^\nu} \Vert_2^2.
\end{aligned}
\]
Consequently, 
\[
\begin{aligned}
0&< (1-\beta)\bDel(x^\nu)+\frac{\kappa_2\kappa_3}{2}(\gamma^{-1}t_\nu)\Vert{d^\nu} \Vert_2^2\\
&\le (1-\beta)\xi+\frac{\kappa_2\kappa_3}{2}(\gamma^{-1}t_\nu)\Vert{d^\nu} \Vert_2^2.
\end{aligned}
\]
Taking the limit over $\nu\in J$ in this inequality gives the contradiction $0\le (1-\beta)\xi <0$.

\end{document}